\documentclass[12pt]{amsart}
\usepackage[top=30truemm,bottom=30truemm,left=25truemm,right=25truemm]{geometry}
\usepackage{txfonts}
\usepackage{mathrsfs}

\usepackage{color}
\usepackage{bm}
\usepackage{amsfonts,amssymb}
\usepackage{dsfont}
\usepackage{extarrows}
\usepackage{amsmath}
\usepackage{mathrsfs}
\usepackage{enumerate}
\usepackage{amscd}
\usepackage[all]{xy}
\usepackage[hyperfootnotes=true]{hyperref}
\theoremstyle{plain} 
\newtheorem{theorem}{\indent\bf Theorem}[section]
\newtheorem{lemma}[theorem]{\indent\bf Lemma}

\theoremstyle{definition} 
\newtheorem{definition}[theorem]{\indent\bf Definition}

\newtheorem{question}[theorem]{\indent\bf Question}
\newtheorem{thm}{Theorem}[section]

\theoremstyle{definition}

\theoremstyle{remark}
\newtheorem{rem}{Remark}[section]

\newcommand{\be}{\begin{equation}}
\newcommand{\ee}{\end{equation}}
\newcommand{\bea}{\begin{eqnarray}}
\newcommand{\eea}{\end{eqnarray}}
\newcommand{\ben}{\begin{eqnarray*}}
	\newcommand{\een}{\end{eqnarray*}}
\newcommand{\bt}{\begin{split}}
	\newcommand{\et}{\end{split}}
\newcommand{\bet}{\begin{equation}}

\begin{document}
	\title[Extensions of K\"ahler currents]{On the extensions of  K\" ahler currents  on compact K\"{a}hler manifolds}
	\author{Zhiwei Wang}
	\address{ Zhiwei Wang:  Laboratory of Mathematics and Complex Systems (Ministry of Education)\\School of Mathematical Sciences\\Beijing Normal University\\Beijing \\ 100875\\P. R. China}
	\email{zhiwei@bnu.edu.cn}
	\author{Xiangyu Zhou}
	\address{Xiangyu Zhou: Institute of Mathematics\\Academy of Mathematics and Systems Sciences\\and Hua Loo-Ken Key Laboratory of Mathematics\\Chinese Academy of Sciences\\Beijing\\100190\\P. R. China}
	\address{School of
		Mathematical Sciences, University of Chinese Academy of Sciences,
		Beijing 100049, China}
	\email{xyzhou@math.ac.cn}
	
	\begin{abstract}
Let $(X,\omega)$ be a compact K\"{a}hler manifold with a K\"{a}hler form $\omega$ of complex dimension $n$, and $V\subset X$ is a compact complex submanifold of positive dimension $k<n$. 
Suppose that $V$ can be embedded in $X$ as a zero section of  a holomorphic vector bundle or rank $n-k$ over $V$.  
Let $\varphi$ be a  strictly $\omega|_V$-psh function  on $V$. 
In this paper, we  prove that there is a strictly $\omega$-psh function $\Phi$  on $X$, such that $\Phi|_V=\varphi$.  
This result  gives a partial answer to an open problem raised by  Collins-Tosatti and Dinew-Guedj-Zeriahi, for the case of K\"{a}hler currents. 
We also discuss possible extensions of K\"ahler currents in a big class.
%

	\end{abstract}

\maketitle
	
\tableofcontents
\section{Introduction}

In this paper, we study the following type extension problem mentioned in Collins-Tosatti  \cite{CT14} and Dinew-Guedj-Zeriahi   \cite[Question 37]{DGZ16}.
\begin{question}\label{DGZ-que}
Let $(X,\omega)$ be a compact K\"{a}hler manifold and $V\subset X$ a complex submanifold.  Is it true that
\begin{align*}
Psh(V,\omega|_V)=Psh(X,\omega)|_V?
\end{align*}
\end{question}

About Question \ref{DGZ-que}  Schumacher \cite{Sch98} proved that if $\omega$ is rational, then any smooth K\"{a}hler form on $V$ in the class $[\omega|_V]$ extends to a smooth K\"{a}hler form on $X$ in the class $[\omega]$.

Coman-Guedj-Zeriahi proved in \cite{CGZ13} that under the same rationality assumption, every $\omega|_V$-psh function $\varphi$ on a closed analytic subvariety $V\subset X$  extends to $\omega$-psh function on $X$.  In the same paper, they also proved that if $(X,\omega)$ is compact K\"{a}hler, any smooth strictly $\omega|_V$-psh function can be extended to a smooth strictly $\omega$-psh function on $X$.

Collins-Tosatti  \cite{CT14} got rid of rationality assumptions in the case of extension of K\"{a}hler currents with analytic singularities from a submanifold. More precisely, they proved that and strictly $\omega|_V$-psh function with analytic singularities extends to a strictly $\omega$-psh function  on $X$.  They  also proved a similar extension result in  \cite{CT15}, which can be used to establish the equality of the non-K\"{a}hler locus  and the null locus of a big and nef $(1,1)$-class on a compact complex  manifold in the Fujiki class $\mathcal{C}$ (i.e. a compact complex manifold bimeromorphic to a compact K\"{a}hler manifold). Recently, this equality was used in  \cite{CHP16,HP16}  to study the minimal models of compact K\"{a}hler 3-folds.

It is mentioned by Dinew-Guedj-Zeriahi \cite{DGZ16} that the general case of  Question \ref{DGZ-que}  is up to now largely open.

In this paper, we answer  Question \ref{DGZ-que} in the case of K\"{a}hler currents partially, by establishing the following
\begin{thm}[Main Theorem]\label{main theorem}
Given a compact K\"{a}hler manifold $(X,\omega)$ of complex dimension $n$, a compact complex submanifold $V\subset X$ of positive dimension $k<n$. 
Suppose that $V$ can be embedded in $X$ as a zero section of a holomorphic vector bundle of rank $n-k$ over $V$. 
Let $\varphi$ be a  strictly $\omega|_V$-psh function. 
Then there is a strictly  $\omega$-psh function $\Phi$  on $X$,  such that $\Phi|_V=\varphi$. 
In other words, any K\"{a}hler  current $T$   in the class $[\omega|_V]$  is a restriction of  a K\"{a}hler   current $\widetilde{T}\in [\omega]$ to $V$.
\end{thm}
We say that $V$ can be embedded in $X$ as a zero section of a holomorphic vector of rank $n-k$ over $V$, if there is a holomorphic Hermitian vector bundle $(E,h)\rightarrow V$ of rank $n-k$, and there is a biholomorphic map $F:E_\delta:=\{(p,v)\in E:|v|_h<\delta\}\rightarrow U$ with $U$ an open neighborhood of $V$ in X, for some small $\delta>0$, such that $F|_V:V:=\{(p,0)\in E\}\rightarrow X$ gives the embedding of $V$ in $X$. 

As a direct consequence, we have the following
\begin{thm}
Given a compact K\"{a}hler manifold $(X,\omega)$ of complex dimension $n$, a compact complex submanifold $V\subset X$ of positive dimension $k<n$. 
Suppose that $V$ can be embedded in $X$ as a zero section of a holomorphic vector bundle of rank $n-k$ over $V$. 
Let $\varphi$ be a  strictly $\omega|_V$-psh function. 
Then for any $0<\varepsilon\ll  1$, there is a $(1+\varepsilon) \omega$-psh function $\Phi$ on $X$, such that $\Phi|_V=\varphi$.
\end{thm}

        Let us give a sketch  the proof of Theorem \ref{main theorem}. For  strictly $\omega|_V$-psh function $\varphi$, we apply Demailly's regularization theorem to get a sequence of decreasing smooth strictly $\omega|_V$-psh functions $\varphi_m$ on $V$. Then we extend uniformly   these  decreasing  $\varphi_m$ on $V$ to decreasing continuous strictly $\omega$-psh functions $\Phi_m$ on $X$ with a uniform  estimate of positivity.  Finally,  from the standard quasi-psh function theory, we know that the limit $\Phi$ of $\Phi_m$  is  a strictly  $\omega$-psh function on $X$, and automatically we  have  that the restriction of $\Phi$ to $V$ is just $\varphi$.

We also discuss the problem of extend K\"ahler currents in a big class on compact K\"ahler manifolds in Section \ref{big}.
\begin{thm}
Given a compact K\"{a}hler manifold $(X,\omega)$ of complex dimension $n$, a compact complex submanifold $V\subset X$ of positive dimension $k<n$. 
Suppose that $V$ can be embedded in $X$ as a zero section of a holomorphic vector bundle of rank $n-k$ over $V$.  Let $\alpha\in H^{1,1}(X,\mathbb R)$ be a big class and $E_{nK}(\alpha)\subset V$. 
Then any  K\"ahler current in  $\alpha|_V$ is the restriction of a K\"ahler current in $\alpha$.
\end{thm}



\subsection*{Acknowledgement}
The first author  would like to thank Prof. Jean-Pierre Demailly for helpful discussion.
Part of the paper was finshed during the the first author's visit of Institut Fourier, and  he would like to thank the hosiptality of the institute. 
The first author was partially supported by the Beijing Natural Science Foundation (1202012, Z190003) and by the NSFC grant (11701031). The second author was partially supported by the NSFC grant (11688101).
	
	
	\section{Facts of quasi-plurisubharmonic functions}
Throughout this paper, we assume that $(X,\omega)$ is a compact K\"{a}hler manifold with K\"{a}hler metric $\omega$ and  $V\subset X$ be a closed complex submanifold of positive dimension,
then  $\omega$ restricts to a K\"{a}hler metric $\omega|_V$ on $V$.
An upper semi-continuous function $\varphi$ is said to be a  quasi-plurisubharmonic function, if $\varphi$ can be written  as $\varphi=\psi+h$ locally, where $\psi$ is a plurisubharmonic function and $h$ is a smooth function.
In particular, a quasi-psh function is in $L^1_{loc}$.
\begin{definition} A quasi-psh function $\varphi$ on $X$ is said to be \begin{itemize}
\item $\omega$-plurisubharmonic, if $\omega+i\partial\bar{\partial}\varphi\geq 0$ in the sense of currents.
\item strictly $\omega$-plurisubharmonic, if $\omega+i\partial\bar{\partial}\varphi>\varepsilon\omega$ in the sense of currents for some $\varepsilon>0$.
\end{itemize}
We denote by  (SPsh$(X,\omega)$) Psh$(X,\omega)$ the set of all the (strictly) $\omega$-psh functions on $X$.
\end{definition}
	
Similar with psh functions, quasi-psh functions also share the following property.
\begin{lemma}[\cite{GZ05}]\label{limit}
Let $\{\varphi_j\}$ be a decreasing family of  non-positive $\omega$-psh functions on a compact K\"ahler manifold $(X,\omega)$, then either $\varphi_j\rightarrow -\infty$ locally uniformly, or $\varphi:=\lim\limits_j\varphi_j\in$ Psh$(X,\omega)$.
\end{lemma}
	
%


Let $T$ be a closed positive  $(1,1)$-current on $X$.
It is proved in \cite{Dem12} that there is a smooth $(1,1)$-form $\alpha$  and a quasi-psh function  $\varphi$ such that $T=\alpha+i\partial\bar{\partial}\varphi$.  Suppose that $\varphi$ is not identically $-\infty$ on $V$, then the   restriction of  $T$ to $V$ is well-defined and $T|_V=\alpha|_V+i\partial\bar{\partial}(\varphi|_V)$, where $\varphi|_V=\varphi\circ \imath$ and  $\imath$ is the holomorphic inclusion map $\imath:V\hookrightarrow X$.

The Lelong number of a quasi-psh function $\varphi$ at a point $x\in X$ is defined as
\begin{align*}
\nu(\varphi,x):=\liminf_{z\rightarrow x}\frac{\varphi(z)}{\log|z-x|}.
\end{align*}
It can also be characterized as  following
\begin{align*}
\nu(\varphi, x)=\sup\{\gamma\geq0:\varphi(z)\leq \gamma\log|z-x|+O(1)   \mbox{~~~at~~~~}  x\}.
\end{align*}

For any quasi-psh function $\varphi$ on $D$, the multiplier ideal sheaf $\mathcal{I}(\varphi)$ of $\varphi$ is defined to be the sheaf of germs of holomorphic functions $f$ such that $|f|^2e^{-2\varphi}$ is locally integrable.

\begin{lemma}[\cite{Sk72}]\label{Sk lemma}
Let $\varphi$ be a quasi-psh function on an open set $\Omega\subset \mathbb{C}^n$ and let $x\in \Omega$.
\begin{itemize}
\item [(a)] If $\nu(\varphi,x)<1$, then $\mathcal{I}(\varphi)_x=\mathcal{O}_{\Omega,x}$.
\item [(b)] If $\nu(\varphi,x)\geq n+s$ for some integer $s\geq 0$, then $\mathcal{I}(\varphi)_x\subset \mathfrak{m}_{\Omega,x}^{s+1}$, where $\mathfrak{m}_{\Omega,x}$ is the maximal ideal of $\mathcal{O}_{\Omega,x}$.
\end{itemize}
\end{lemma}

	The following celebrated regularization theorem is due to Demailly.
	\begin{thm}[c.f. \cite{Dem12,DemA}]\label{regularization}
		Let $\varphi$ be a quasi-psh function on a compact Hermitian manifold $(X,\omega)$ such that $\frac{i}{\pi}\partial\bar{\partial}\varphi\geq \gamma$ for some continuous $(1,1)$-form $\gamma$. Then there is a sequence of quasi-psh functions $\varphi_m$ with logarithmic singularities and a decreasing sequence $\varepsilon_m>0$ converging to $0$ such that
		\begin{itemize}
			\item [(a)] $\varphi(x)<\varphi_m(x)\leq \sup_{|\zeta-x|<r}\varphi(\zeta)+C\Big( \frac{|\log r|}{m}+r+\varepsilon_m   \Big)$ with respect to coordinate open sets covering $X$. In particular, $\varphi_m$ converges to $\varphi$ pointwise and in $L^1(X)$ and
			\item [(b)] $\nu(\varphi,x)-\frac{n}{m}\leq \nu(\varphi_m,x)\leq \nu(\varphi,x)$ for every $x\in X$;
			\item [(c)] $\frac{i}{\pi}\partial\bar{\partial}\varphi_m\geq \gamma-\varepsilon_m\omega$.
		\end{itemize}
Furthermore, for any multiplicative subsequence $m_k$, one can arrange that $\varphi_{m_k}$ is a non-increasing sequence of potentials.
	\end{thm}
If the given quasi-psh function has zero Lelong number everywhere,  from Lemma \ref{Sk lemma}, we see that $\varphi_m$ are smooth for all $m$.
	
Applying  Theorem \ref{regularization} to every $\varphi_j:=\max\{\varphi,-j\}+\frac{1}{j}$ with $j\in \mathbb{Z}^+$, and by a diagonal argument, we can  obtain the following
\begin{lemma}[\cite{Blocki-Kolodziej07}]\label{key lemma}
Let $\varphi$ be a quasi-psh function on a compact Hermitian  manifold $(X,\omega)$,  such that $\omega+i\partial\bar{\partial}\varphi\geq \varepsilon\omega$.
Then there is a  sequence of smooth functions $\varphi_m$ and a decreasing sequence $\varepsilon_m>0$ converging to $0$, satisfying the following
\begin{itemize}
\item [(a)] $\varphi_m\searrow \varphi$;
\item [(b)]$\omega+i\partial\bar{\partial}\varphi_m\geq (\varepsilon-\varepsilon_m)\omega$.
\end{itemize}
\end{lemma}
\begin{rem}We note that similar results was obtained by B\l ocki-Kolodziej \cite{Blocki-Kolodziej07}  using more elementary methods.
\end{rem}

Finally,  we state the following lemma   due to  Demailly-P\u{a}un for later use.
\begin{lemma}[c.f. \cite{DP04}]\label{reference function}
There exists a function $F:X\rightarrow [-\infty, +\infty)$ which is smooth on $X\setminus V$, with logarithmic singularities along $V$, and such that $\omega+i\partial\bar{\partial}F\geq \varepsilon \omega$ is  a K\"{a}hler current on $X$.
By subtracting a large constant, we can make that $F<0$ on $X$.
\end{lemma}
\section{Extension of  strictly  $\omega|_V$-psh functions: proof of the Main Theorem }\label{distance function}
In this section, we will present the proof of Theorem \ref{main theorem}, which divides into three steps.
	
Let $T=\omega|_V+i\partial\bar{\partial}\varphi\geq \varepsilon\omega|_V$ be the given K\"{a}hler current  in the K\"{a}hler class $[\omega|_V]$, where $\varphi$ is a strictly $\omega|_V$-psh function.
By subtracting a large constant, we may assume that $\sup_V \varphi<-C$ for some positive constant $C$.

By  Lemma \ref{key lemma}, we have that there is a sequence of non-increasing  smooth strictly $\omega|_V$-psh functions $\varphi_{m}$ on $V$,
and a decreasing sequence of positive numbers $\varepsilon_m$ such that as $m\rightarrow \infty$
\begin{itemize}
\item $\varepsilon_m\rightarrow 0$;
\item $\varphi_{m} \searrow \varphi$;
\item $\omega|_V+i\partial\bar{\partial}\varphi_{m}    \geq (\varepsilon-\varepsilon_m)\omega|_V$.
	\end{itemize}
Since $\varphi<-C$ on $V$,  from (a) in Lemma  \ref{key lemma}, one can see that for $m$ large, $\varphi_{m}<-\frac{C}{2}$. By choosing  a subsequence, we assume that for any $m\in \mathbb{N}$, $\varphi_{m}<-\frac{C}{2}$ and $\omega|_V+i\partial\bar{\partial}\varphi_m> \frac{\varepsilon}{2}\omega|_V$.
	
We say a smooth strictly $\omega|_V$-psh function $\phi$ on $V$ satisfies \textbf{assumption $\bigstar_{\varepsilon, C}$}, if $\omega|_V+i\partial\bar{\partial}\phi>\frac{\varepsilon}{2}\omega|_V$ and $\phi<-\frac{C}{2}$.
	
Note that for all $m\in \mathbb{N}^+$, $\varphi_m$ satisfy  \textbf{assumption $\bigstar_{\varepsilon, C}$}.  In the following, we will extend all the $\varphi_m$ simultaneously to non-increasing  strictly $\omega$-psh functions on   the ambient manifold $X$.
	
\subsection*{Step 1: Local uniform extensions of $\varphi_m$ for all $m$ }
Let $\phi$ be a function satisfying  \textbf{assumption $\bigstar_{\varepsilon, C}$}. Since $V$ is embedded in $X$ as a zero section of a holomorphic vector bundle, say $E$, there is a open subset $U$ of $V$ in $X$, biholomorphic to a neighbourhood of the zero section $V$ in $E$. 
Let $\pi:E\rightarrow V$ be the natural holomorphic projection of the holomorphic vector bundle $E$, and $h$ be a hermitian metric on $E$. There is an induced " square distance function" $\rho$ defined on $U$ as
\begin{align*}
\rho(v_p):=h(v_p,v_p)
\end{align*}
for any $v_p\in E_p$.
From the holomorphic vector bundle structure, for any $p\in V$,   there is a  holomorphic  local coordinate chart $(U_p,(z_1,\cdots,z_{n-k},z_{n-k+1},\cdots,z_n))$ in $U$ centered at $p$, such that $(z_{n-k+1},\cdots, z_n)$ serve as holomorphic coordinates on $V$ near $p$,  and $(z_1,\cdots, z_{n-k})$ a holomorphic coordinates of the fibers of $E$ with respect to a holomorphic basis $e_1,\cdots, e_{n-k}$ of $E$ near $p$.  
Write $z=(z_1,\cdots, z_{n-k})$ and $z'=(z_{n-k+1},\cdots,z_n)$, then we have $\pi(z,z')=(0,z')$ and  $\rho=\sum_{i,j=1,\cdots, n-k}h_{i\bar j}(z')z_i\bar z_j$.

On $U_p$,  we define 	
\begin{align*}
\bar{\phi}(z,z'):=(\phi\circ \pi)(z,z')+A\rho(z,z')=\phi(z')+A\sum_{i,j=1,\cdots, n-k}h_{i\bar j}(z')z_i\bar z_j
\end{align*}
where $A$ is a positive constant to be determined later. 
	

We have the following computations:
\begin{align*}
\omega+i\partial\bar{\partial}\bar{\phi}&=\sum_{i,j=1,\cdots,n-k}g_{i\bar{j}}(z,z')dz_i\wedge d\bar{z}_j +Ai\sum_{i=1}^{n-k}h_{i\bar j}(z')dz_i\wedge d\bar{z}_i\\
&+\sum_{i=1,\cdots,n-k;j=n-k+1,\cdots,n}\Big(g_{i\bar{j}}(z,z')+A\sum_{1\leq l\leq n-k}\frac{\partial h_{i\bar k}(z')}{\partial \bar z_j}\bar z_l\Big)dz_i\wedge d\bar{z}_j \\
&+\sum_{i=n-k+1,\cdots,n;j=1,\cdots,n-k}\Big(g_{i\bar{j}}(z,z')+A\sum_{1\leq l\leq n-k}\frac{\partial h_{l\bar{j}}(z')}{\partial z_i}z_l\Big)dz_i\wedge d\bar{z}_j\\
&+\sum_{i,j=n-k+1,\cdots,n}g_{i\bar{j}}(z,z')dz_i\wedge d\bar{z}_j+i\partial\bar{\partial}\phi(z')\\
&=I+II+III+IV.
\end{align*}
In matrix form, we write $\omega+i\partial\bar{\partial}\bar{\phi}$ as
\begin{align*}
\Big(	\begin{array}{cc}
I, & II\\
III, & IV
\end{array}\Big).
\end{align*}
It is obvious that $II=III^*$, where $III^*$ is the conjugate transpose of $III$.

One can see that $I$ is positive definite even if $A=0$.
Apply a congruent transformation which is independent of $\phi$, to the above matrix, one can get that
\begin{align*}
\Big(	\begin{array}{cc}
I, & 0\\
0, & IV-(III) I^{-1}(II)
\end{array}\Big).
\end{align*}
Now we have to deal with the term $ IV-(III) I^{-1}(II)$.
Note that
\begin{align*}
IV&=\Big(g_{i\bar{j}}(z,z')+\partial_{z_i}\bar{\partial}_{z_j}\phi(z')\Big)_{i,j=n-k+1,\cdots,n}\\
&=\Big(g_{i\bar{j}}(z,z')-g_{i\bar{j}}(0,z')\Big)_{i,j=n-k+1,\cdots,n}\\
&+\Big(g_{i\bar{j}}(0,z')+i\partial_{z_i}\bar{\partial}_{z_j}\phi(z')\Big)_{i,j=n-k+1,\cdots,n}.
\end{align*}
Considering  the matrix
$$\Big(g_{i\bar{j}}(0,z')+i{\partial_{z_i}\bar{\partial}_{z_j}}\phi(z')\Big)_{i,j=n-k+1,\cdots,n}$$
as an  extension of the  matrix of
$$\omega|_V+i\partial\bar{\partial}\phi(z')$$
from $V\cap U$ to $U$, which,  from assumption,
$$\geq \varepsilon\Big(g_{i\bar{j}}(0,z')\Big)_{i,j=n-k+1,\cdots,n}.
$$
Since on $V$,  $$\Big(g_{i\bar{j}}(0,z')\Big)_{i,j=n-k+1,\cdots,n}$$ is positive definite as a $k\times k$ Hermitian matrix, and  from the smoothness, we have that up to shrinking (note that the shrinking is independent of $\phi$), it is positive definite on $U$.

For the matrix  $$\Big(g_{i\bar{j}}(x,z')-g_{i\bar{j}}(0,z')\Big)_{i,j=n-k+1,\cdots,n},$$from the smoothness of $\omega$, up to shrinking (note that the shrinking is independent of $\phi$), the following holds on $U$:
\begin{align*}
-\frac{\varepsilon}{4}\Big(g_{i\bar{j}}(0,z')\Big)_{i,j=n-k+1,\cdots,n}&\leq \Big(g_{i\bar{j}}(x,z')-g_{i\bar{j}}(0,z')\Big)_{i,j=n-k+1,\cdots,n}\\
&\leq \frac{\varepsilon}{4}\Big(g_{i\bar{j}}(0,z')\Big)_{i,j=n-k+1,\cdots,n}.
\end{align*}

Now on $U$, we get that
\begin{align*}
IV-(III) I^{-1}(II)\geq \frac{3\varepsilon}{4}\Big(g_{i\bar{j}}(0,z')\Big)_{i,j=n-k+1,\cdots,n} -(III) I^{-1}(II).
\end{align*}

As for the matrix 
$$(III) I^{-1}(II),$$
note that $\sum_{1\leq l\leq n-k}\frac{\partial h_{i\bar k}(z')}{\partial \bar z_j}\bar z_l=0$ and $\sum_{1\leq l\leq n-k}\frac{\partial h_{l\bar{j}}(z')}{\partial z_i}z_l=0$ at $p$, and 
 is totally independent of $\phi$, we first choose sufficiently large $A$, and then shrink $U_p$ (independent of $\phi$) if necessary,  such that on $U_p$ $$(III) I^{-1}(II)\leq  \frac{\varepsilon}{4}\Big(g_{i\bar{j}}(0,z')\Big)_{i,j=n-k+1,\cdots,n}.$$
Then we  get an open neighborhood $U_p$  near $p$, which is independent of $\phi$,  such that on $U$, the $k\times k$ matrix $$IV-(III) I^{-1} (II)\geq \frac{\varepsilon}{2}\Big(g_{i\bar{j}}(0,z')\Big)_{i,j=n-k+1,\cdots,n}.$$

Putting all the computations together, we see that
\begin{align*}
\Big(	\begin{array}{cc}
I, & 0\\
0, & IV-(III) I^{-1}(II)
\end{array}\Big)\geq \Big(	\begin{array}{cc}
I, & 0\\
0, & \frac{\varepsilon}{2}\Big(g_{i\bar{j}}(0,z')\Big)_{i,j=n-k+1,\cdots,n}
\end{array}\Big)\geq \varepsilon'\omega
\end{align*}
for some positive constant $\varepsilon'$ (independent of  $\phi$). Moreover, we can shrink $U_p$ (the shrinking is independent of $\phi$), such that $\bar\phi\leq -\frac{C}{4}$.

To emphasis the uniformity, it is worth to point out again that the chosen of the open set $U$, the large constant $A$ and the constant $\varepsilon'$ is independent of $\phi$, as long as $\phi$ satisfies \textbf{assumption $\bigstar_{\varepsilon,C}$}.
We call the above data $ (U,A,\varepsilon',-\frac{C}{4},\bar\phi)$ an \textbf{admissible local extension} of $\phi$.

Since all the $\varphi_m$ satsifies the same \textbf{assumption $\bigstar_{\varepsilon,C}$}, thus near $p$, we can choose a \textbf{uniform admissible local extension $ (U_p,A,\varepsilon',-\frac{C}{4},\bar\varphi_m)$} of $\varphi_m$, for all $m\in\mathbb{N}^+$.
 Since $V$ is compact, one may choose an  open neighborhood $U$ of $V$ in $X$, and universal constants $A>0$ and $\varepsilon'>0$, such that the functions $\widetilde \varphi_m:=\varphi\circ \pi+A\rho$ are defined on $U$, such that $\omega+i\partial\bar\partial \widetilde \varphi_m\geq \varepsilon'\omega$ on $U$ for all $m$. Since $\{\varphi_m\}$ is a non-increasing sequence, one obtain that $\{\widetilde{\varphi}_m\}$ is a non-increasing sequence.

\subsection*{Step 2: Global extensions of $\varphi_m$ for all $m$} Up to shrinking, we may assume that $\widetilde{\varphi}_m$ are defined on the closure of $U$ for all $m\in \mathbb{N}^+$.
Let $F$  be the quasi-psh function in Lemma \ref{reference function}.
Near $\partial U$ (the boundary of $U$), the function $F$ is smooth, and $\sup_{\partial U}\widetilde{\varphi}_{1}=-C''$ for some positive constant $C''>0$.
Now we choose a small positive $\nu$, such that $\inf_{\partial U}(\nu F)>-\frac{C''}{2}$ and   $\omega+i\partial\bar{\partial}\nu F\geq\varepsilon'\omega$.
Thus $\nu F >\widetilde{\varphi}_{1}\geq \widetilde{\varphi}_m$ in a neighborhood of $\partial U$ for all $m\in \mathbb{N}^+$, since $\widetilde{\varphi}_m$ is non-increasing.
Therefore, we can finally define
\begin{align*}
\Phi_m=\left\{
\begin{array}{ll}
\max\{\widetilde{\varphi}_m, \nu F\}, & \hbox{on $U$;} \\
\nu F, & \hbox{on $X\setminus U$,}
\end{array}
\right.
\end{align*}
which is defined on the whole of $X$. It  is easy to check that $\Phi_m$ satisfies the following properties:
\begin{itemize}
\item $\Phi_m$ is non-increasing in  $m$,
\item $\Phi_m\leq 0$ for all $m\in \mathbb{N}^+$,
\item $\omega+i\partial\bar{\partial}\Phi_m\geq \varepsilon'\omega$ for all $m\in \mathbb{N}^+$,
\item $\Phi_m|_V=\varphi_m$ for all $m\in \mathbb{N}^+$.
\end{itemize}
	
\subsection*{Step 3: Taking limit to complete the proof of Theorem \ref{main theorem}}
From above steps, we get a sequence of non-increasing, non-positive strictly $\omega$-psh functions $\Phi_m$ on $X$. Then from Lemma \ref{limit},  we conclude that either $\Phi_m\rightarrow -\infty $ uniformly on $X$, or $\Phi:=\lim\limits_m\Phi_m\in$ Psh$(X,\omega)$.
But $\Phi_m|_V=\varphi_m\searrow \varphi\not\equiv -\infty$, the first case will not appear.
Moreover,   we can see that $\Phi:=\lim\limits_m\Phi_m$ is a strictly $\omega$-psh function on $X$ from the property $\omega+i\partial\bar{\partial}\Phi_m\geq \varepsilon'\omega$ for all $m\in \mathbb{N}^+$, and $\Phi|_V=\lim\limits_m\Phi_m|_V=\lim\limits_m\varphi_m=\varphi$.
It follows that   $(\omega+i\partial\bar{\partial}\Phi)|_V=\omega|_V+i\partial\bar{\partial}\varphi$.
Thus we complete the proof of Theorem \ref{main theorem}.


\begin{rem}
The assumption that $V$ can be embedded in $X$ as a zero section of a holomorphic vector bundle of rank $n-k$ over $V$, is related to the understanding  of the germ, or the neighborhood structure of  a submanifold $V$  in $X$,  which is was studied extensively, e.g. see \cite{Gra62,Gri66,CG79,Hir81, Hwa19}.
\end{rem}

\section{Extensions of K\"ahler currents in a big class}	\label{big}

Let $(X,\omega)$ be a compact K\"ahler manifold of complex dimension $n$, and $V\subset X$ be a complex submanifold of complex dimension $0<k<n$. 
Suppose that $V$ can be embedded in $X$ as a zero section of a holomorphic vector bundle of rank $n-k$ over $V$. 
Let $\alpha\in H^{1,1}(X,\mathbb R)$ be a big class, and $\theta\in \alpha$ be a smooth representative. 
Let $\varphi$ be a quasi-psh function on $V$ such that $\theta|_V+i\partial\bar\partial \varphi\geq \varepsilon\omega|_V$ on $V$, for some $\varepsilon>0$. 
Then by the same technique as in Step 1 and Step 2, we can get an open neighbourhood $U$ of $V$ in $X$, and a non-increasing sequence of smooth functions $\widetilde \varphi_m$ on $U$, such that
\begin{itemize}
\item $\theta+i\partial\bar\partial \widetilde \varphi_m\geq \varepsilon'\omega$ on $U$,
\item $\widetilde{\varphi}_m|_V=\varphi$.
\end{itemize}
In \cite{Bou04}, the non-K\"ahler locus of $\alpha$, is defined as
$$E_{nK}(\alpha):=\cap_{T\in \alpha}E_+(T),$$
where $E_+(T)$ is the set of points of $X$ such that the K\"ahler current $T$ has positive Lelong numbers, and $T$ varies in all the K\"ahler currents in $\alpha$. 
From Siu's semicontinuity of Lelong number upper level sets and strong Noether property, $E_{nK}(\alpha)$ is an analytic subvariety. 
Boucksom \cite[Theorem 3.17]{Bou04} proved that there is a K\" ahler currents $T$ with analytic singularities in $\alpha$, such that $E_+(T)=E_{nK}(\alpha)$. 
We write $T=\theta+i\partial\bar\partial \Upsilon$. Suppose now that $E_{nK}(\alpha)$ is contained in $V$, we set
\begin{align*}
\Phi_m=\left\{
\begin{array}{ll}
\max\{\widetilde{\varphi}_m,\Upsilon+\nu F+C\}, & \hbox{on $U$;} \\
\Upsilon+\nu F+C, & \hbox{on $X\setminus U$,}
\end{array}
\right.
\end{align*}
where $\nu>0$ and $C>0$ is chosen to make  $\inf_{\partial U}\Upsilon+\nu F+C\geq \sup_{\partial U}\widetilde{\varphi}_1$, and $\theta+i\partial\bar\partial(\Upsilon+\nu F+C)\geq\varepsilon''\omega$,  up to shrinking $U$ if necessary.
Therefore, we get a non-increasing sequence of continuous strictly $\theta$-psh functions on $X$, and by the same argument as in Step 3, we conclude that $\Phi:=\lim \Phi_m$ is a desired extension of $\varphi$. We arrived at the following 
\begin{thm}
Let $(X,\omega)$ be a compact K\"ahler manifold of complex dimension $n$, and $V\subset X$ be a complex submanifold of positive dimension. Suppose that $V$ can be embedded in $X$ as a zero section of a holomorphic vector bundle of rank $n-k$ over $V$.  Let $\alpha\in H^{1,1}(X,\mathbb R)$ be a big class and $E_{nK}(\alpha)\subset V$. 
Then any  K\"ahler current in  $\alpha|_V$ is the restriction of a K\"ahler current in $\alpha$.
\end{thm}

\end{document}